\documentclass[11pt]{amsart}
\usepackage{amssymb, amsthm, amscd, amsmath}

\newtheorem{theorem}{Theorem}[section]
\newtheorem{lemma}[theorem]{Lemma}
\newtheorem{prop}[theorem]{Proposition}
\theoremstyle{definition}

\newtheorem{conj}[theorem]{Conjecture}
\newtheorem{remark}[theorem]{Remark}
\newtheorem{example}[theorem]{Example}
 
\def\image{\mathop{\mathrm{Image}}}

\begin{document}

\title{$p$-torsion elements in local cohomology modules II}

\author{Anurag K. Singh}
 
\address{Department of Mathematics, University of Utah, 155 S. 1400 E., Salt
Lake City, \quad UT 84112, USA \newline E-mail: {\tt singh@math.utah.edu}}

\thanks{February 14, 2001 \\
This manuscript is based on work supported in part by the National Science
Foundation under Grant No. DMS 0070268}

\maketitle

\section{Introduction}

Gennady Lyubeznik conjectured that if $R$ is a regular ring and ${\mathfrak a}$
is an ideal of $R$, then the local cohomology modules $H_{\mathfrak a}^i(R)$
have only finitely many associated prime ideals, \cite[Remark 3.7 (iii)]{Ly1}.
While this conjecture remains open in this generality, several results are now
available: if the regular ring $R$ contains a field of prime characteristic $p
> 0$, Huneke and Sharp showed in \cite{HS} that the set of associated prime
ideals of $H_{\mathfrak a}^i(R)$ is finite. If $R$ is a regular local ring
containing a field of characteristic zero, Lyubeznik showed that $H_{\mathfrak
a}^i(R)$ has only finitely many associated prime ideals, see \cite{Ly1} and
also \cite{Ly2, Ly3}. Recently Lyubeznik has also proved this result for
unramified regular local rings of mixed characteristic, \cite{Ly4}.

In \cite{Hu} Craig Huneke first raised the following question: for a Noetherian
ring $R$, an ideal ${\mathfrak a} \subset R$, and a finitely generated
$R$-module $M$, is the number of associated primes ideals of $H_{\mathfrak
a}^i(M)$ always finite? For some of the work on this problem, we refer the
reader to the papers \cite{BL, BRS, He} in addition to those mentioned above.
In \cite{Si} we constructed an example of a hypersurface $R$ for which a local
cohomology module $H_{\mathfrak a}^3(R)$ has $p$-torsion elements for every
prime integer $p$, and consequently has infinitely many associated prime
ideals. Since this is the only known source of infinitely many associated prime
ideals so far, it is worthwhile to investigate whether similar techniques may
yield an example of a regular ring $R$ for which a local cohomology module
$H_{\mathfrak a}^i(R)$ has $p$-torsion elements for every prime integer $p$.
This leads to some very intriguing questions as we shall see in this paper. Our
results thus far support Lyubeznik's conjecture that local cohomology modules
of all regular rings have only finitely many associated prime ideals. 

Let $R$ be a polynomial ring over the integers and $F_i,G_i$ be elements of $R$
for which
$$
F_1G_1 + F_2G_2 + \dots + F_nG_n = 0.
$$
Consider the ideal ${\mathfrak a} = (G_1, \dots, G_n)R$ and the local 
cohomology module 
$$
H_{\mathfrak a}^n(R) = \varinjlim \frac{R}{(G_1^k, \dots, G_n^k)R}
$$ 
where the maps in the direct limit system are induced by multiplication by the 
element $G_1 \cdots G_n$. For a prime integer $p$ and prime power $q=p^e$, let
$$
\lambda_q = \frac{(F_1G_1)^q + \dots + (F_nG_n)^q}{p}.
$$ 
Note that has $\lambda_q$ has integer coefficients, i.e., that 
$\lambda_q \in R$. Consider
$$
\eta_q = [\lambda_q + (G_1^q, \dots, G_n^q)R] \ \in \ H_{\mathfrak a}^n(R)
= \varinjlim \frac{R}{(G_1^k, \dots, G_n^k)R}. 
$$
It is immediately seen that $p\eta_q=0$ and so if $\eta_ q$ is a nonzero 
element of $H_{\mathfrak a}^n(R)$, then it must be a $p$-torsion element.
Hence if the local cohomology module $H_{\mathfrak a}^n(R)$ has only finitely
many associated prime ideals then, for all but finitely many prime integers
$p$, the elements $\eta_q$ as constructed above must be zero i.e., for some
$k \in \mathbb N$, we have 
$$
\lambda_q (G_1 \cdots G_n)^{k} \in (G_1^{q+k}, \dots, G_n^{q+k})R.
$$

This motivates the following conjecture:
\begin{conj}
Let $R$ be a polynomial ring over the integers and let $F_i,G_i\in R$ such that 
$$
F_1G_1 + \dots + F_nG_n = 0.
$$
Then for every prime integer $p$ and prime power $q=p^e$, there exists 
$k \in \mathbb N$ such that 
$$
\frac{(F_1G_1)^q + \dots + (F_nG_n)^q}{p} (G_1 \cdots G_n)^{k} \in
(G_1^{q+k}, \dots, G_n^{q+k})R.
$$
\label{conj}
\end{conj}

We shall say that Conjecture~\ref{conj} holds for $G_1,\dots,G_n \in R$, if it
holds for all relations $\displaystyle{\sum_{i=0}^n F_i G_i =0}$ with
$F_1,\dots,F_n \in R$. 

\begin{remark} 
The above conjecture is easily established if $n=2$ since in this case we have 
$F_1G_1 + F_2G_2 = 0$ and so 
$$
\frac{(F_1G_1)^q + (F_2G_2)^q}{p} = \frac{(F_1G_1)^q + (-F_1G_1)^q}{p} =
\begin{cases}
(F_1G_1)^q & \text{if } \ p=2, \\
0 & \text{if } \ p \neq 2,
\end{cases}
$$
which is an element of $(G_1^q, G_2^q)R$.
\end{remark}

\begin{example}{\bf The hypersurface example.} 
Conjecture~\ref{conj} is false if the condition that $R$ is a polynomial ring
over the integers is replaced by the weaker condition that the ring $R$ is a
hypersurface over the integers. To see this, let
$$
R={\mathbb Z}[U,V,W,X,Y,Z]/(UX+VY+WZ)
$$
and ${\mathfrak a}$ denote the ideal $(x,y,z)R$. (We use lowercase letters here 
to denote the images of the corresponding variables.) Consider the relation 
$$
ux+vy+wz=0
$$
where, in the notation of Conjecture~\ref{conj}, $F_1=u, \ F_2=v, \ F_3=w$ and 
$G_1=x, \ G_2=y, \ G_3=z$. By means of a multi-grading argument, it is 
established in \cite[\S 4]{Si} that 
$$
\frac{(ux)^p + (vy)^p + (wz)^p}{p}(xyz)^k \notin 
(x^{p+k}, \ y^{p+k}, \ z^{p+k})R \ \text { for all } \ k \in {\mathbb N}.
$$
Consequently Conjecture~\ref{conj} does not hold here even for the choice of the
prime power $q=p$. This is, of course, the example from \cite{Si} of a local
cohomology module with infinitely many associated prime ideals: more precisely,
for every prime integer $p$, 
$$
\left[ \frac{(ux)^p + (vy)^p + (wz)^p}{p} + \ (x^p, y^p, z^p)R \right] \ 
\in \ H^3_{(x,y,z)}(R)
$$
is a nonzero $p$-torsion element. 
\end{example}

\begin{conj}
We record another formulation of Conjecture~\ref{conj}. Let $F_i,G_i$ be
elements of a polynomial ring $R$ over the integers such that 
$$
F_1G_1 + \dots + F_nG_n = 0.
$$
Let ${\mathfrak a}$ denote the ideal $(G_1, \dots, G_n)$ of $R$. For an 
arbitrary prime integer $p$ and prime power $q=p^e$, we have
$$
(F_1G_1)^q + \dots + (F_nG_n)^q \equiv 0 \mod p
$$ 
which is a relation on the elements $\overline{G}_1^q, \dots, \overline{G}_n^q 
\in R/pR$, where $\overline{\phantom{G}}$ denotes the image of an element of
$R$ in the ring $R/pR$. This relation may be viewed as an element $\mu_q \in 
H_{\mathfrak a}^{n-1}(R/pR)$. Conjecture~\ref{conj} is equivalent to the
conjecture that this element $\mu_q$ is in the image of the natural 
homomorphism $H_{\mathfrak a}^{n-1}(R) \to H_{\mathfrak a}^{n-1}(R/pR)$.

To see the equivalence of these conjectures, suppose 
$$
\mu_q \in \image\Big( H_{\mathfrak a}^{n-1}(R) 
 \to H_{\mathfrak a}^{n-1}(R/pR) \Big).
$$
Then the relation $(\overline{F}_1^q, \dots, \overline{F}_n^q)$ on the elements 
$\overline{G}_1^q, \dots, \overline{G}_n^q \in R/pR$ lifts to a relation in 
$H_{\mathfrak a}^{n-1}(R)$, i.e., there exists an integer $k$ and elements 
$\alpha_i \in R$ such that 
\begin{align*}
& \alpha_1 G_1^{q+k}+ \dots + \alpha_n G_n^{q+k}=0 \qquad \text{ and} \\
& \alpha_i \equiv F_i^q (G_1 \cdots G_{i-1} G_{i+1} \cdots G_n)^k \mod p
\quad \text{ for all } \quad 1 \le i \le n.
\end{align*}
Hence we have
\begin{multline*}
\Big( (F_1G_1)^q + \dots + (F_nG_n)^q \Big)(G_1 \dots G_n)^k \\
= (F_1^q G_2^k \cdots G_n^k -\alpha_1)G_1^{q+k} + \dots + 
        (F_n^q G_1^k \cdots G_{n-1}^k -\alpha_n)G_n^{q+k} \\
\in \Big( p G_1^{q+k}, \dots, p G_n^{q+k} \Big)R,
\end{multline*}
and so 
$\displaystyle{
\frac{(F_1G_1)^q + \dots + (F_nG_n)^q}{p} (G_1 \cdots G_n)^{k} \in
\Big( G_1^{q+k}, \dots, G_n^{q+k} \Big)R}$.

The proof of the converse is similar.
\end{conj}

\begin{remark}
We next mention a conjecture due to Mel Hochster. While this was shown to be
false in \cite{Si}, our entire study of $p$-torsion elements originates from
this conjecture.

Consider the polynomial ring over the integers $R={\mathbb Z}[u,v,w,x,y,z]$
where ${\mathfrak a}$ is the ideal generated by the size two minors of the
matrix 
$$
M = \begin{pmatrix} 
u & v & w \\ 
x & y & z \\ 
\end{pmatrix}, 
$$
i.e., ${\mathfrak a} = (\Delta_1, \Delta_2, \Delta_3)R$ where $\Delta_1 =
vz-wy$, $\Delta_2 = wx-uz$, and $\Delta_3 = uy-vx$. In the ring $R$ we have
$u\Delta_1 + v\Delta_2 + w\Delta_3 = 0$. The relation
$$
(u\Delta_1)^q + (v\Delta_2)^q + (w\Delta_3)^q \equiv 0 \mod p
$$ 
may be viewed as an element $\mu_q \in H_{\mathfrak a}^2(R/pR)$. Hochster 
conjectured that for every prime integer there exists a choice of $q=p^e$ such 
that
$$
\mu_q \notin \image \Big(H_{\mathfrak a}^2(R) \to H_{\mathfrak a}^2(R/pR)\Big),
$$
and consequently that the image of $\mu_q$ in $H_{\mathfrak a}^3(R)$ is a 
nonzero $p$-torsion element of $H_{\mathfrak a}^3(R)$. In \cite{Si} we
constructed an equational identity which provided us with an element of 
$H_{\mathfrak a}^2(R)$ that maps to $\mu_q \in H_{\mathfrak a}^2(R/pR)$.
While we refer the reader to \cite{Si} for the details of the construction, we 
would like to provide a brief sketch. 

Consider the following equational identity:
\begin{align*}
&\Delta_1^{2k+1} u^{k+1} \sum_{n=0}^k \binom{k}{n} x^n \sum_{i=0}^n (-1)^i
\binom{k+i}{k} \binom{k+n-i}{k} w^i v^{n-i} \Delta_2^{k-i}\Delta_3^{k-n+i} \\
+&\Delta_2^{2k+1} v^{k+1} \sum_{n=0}^k \binom{k}{n} y^n \sum_{i=0}^n (-1)^i
\binom{k+i}{k} \binom{k+n-i}{k} u^i w^{n-i} \Delta_3^{k-i}\Delta_1^{k-n+i} \\
+&\Delta_3^{2k+1} w^{k+1} \sum_{n=0}^k \binom{k}{n} z^n \sum_{i=0}^n (-1)^i
\binom{k+i}{k} \binom{k+n-i}{k} v^i u^{n-i} \Delta_1^{k-i}\Delta_2^{k-n+i} \\
&=0.
\end{align*}
As it is a relation on the elements $\Delta_i^{2k+1}$ this identity gives us, 
for every $k \in \mathbb N$, an element $\gamma_k \in H_{\mathfrak a}^2(R)$. 
Using $k=q-1$ and examining the binomial coefficients above $\mod p$, we obtain 
$$
\Big((u\Delta_1)^q+(v\Delta_2)^q+(w\Delta_3)^q\Big)
(\Delta_1 \Delta_2 \Delta_3)^{q-1} \equiv 0 \mod p.
$$
Consequently $\gamma_{q-1} \mapsto \mu_q$ under the natural homomorphism
$$
H_{\mathfrak a}^2(R) \to H_{\mathfrak a}^2(R/pR).
$$
\end{remark}

\begin{prop}
If Conjecture~\ref{conj} is true for the relations 
$\displaystyle{\sum_{i=0}^n E_i G_i =0}$ and 
$\displaystyle{\sum_{i=0}^n F_i G_i =0}$ where
$E_i, F_i, G_i \in R={\mathbb Z}[X_1,\dots,X_m]$, then it is also true for 
the relation $\displaystyle{\sum_{i=0}^n (sE_i+tF_i) G_i =0}$ in 
$S={\mathbb Z}[X_1,\dots,X_m,s,t]$. More precisely, if for a prime 
power $q=p^e$, there exists $k_1, k_2 \in \mathbb N$ such that 
\begin{multline*}
\left[ \sum_{i=0}^n \frac{(E_iG_i)^q}{p} \right] (G_1 \cdots G_n)^{k_1} \in
   \Big( G_1^{q+k_1}, \dots, G_n^{q+k_1} \Big) R \ \text{ and } \\ 
\left[ \sum_{i=0}^n \frac{(F_iG_i)^q}{p} \right] (G_1 \cdots G_n)^{k_2} \in
  \Big( G_1^{q+k_2}, \dots, G_n^{q+k_2} \Big) R, \ \text{ then } \\ 
\left[ \sum_{i=0}^n \frac{(sE_i+tF_i)^q G_i^q}{p} \right] (G_1 \cdots G_n)^{k} 
 \in \Big( G_1^{q+k}, \dots, G_n^{q+k} \Big) S 
\end{multline*}
for $k = \max\{k_1,k_2\}$.
\label{sum}
\end{prop}

We leave the proof as an elementary exercise. 

\section{A special case of the conjecture}

In Theorem~\ref{reg} below we prove what is perhaps the first interesting case
of Conjecture~\ref{conj}. 

\begin{theorem}
Let $R$ be a polynomial ring over the integers and $F_i, G_i$ be elements of $R$
such that $F_1, F_2, F_3$ form a regular sequence in $R$ and
$$
F_1G_1 + F_2G_2 + F_3G_3 = 0.
$$
Let $q=p^e$ be a prime power. Then for $k=q-1$ we have
$$
\frac{(F_1G_1)^q + (F_2G_2)^q + (F_3G_3)^q}{p} (G_1 G_2 G_3)^k \ \in \
(G_1^{q+k}, G_2^{q+k}, G_3^{q+k})R.
$$
\label{reg}
\end{theorem}

\begin{remark}
In Hochster's conjecture discussed earlier, the relation under consideration is
$u\Delta_1 + v\Delta_2 + w\Delta_3 = 0$ and the elements $u,v,w$ certainly
form a regular sequence in the polynomial ring $R={\mathbb Z}[u,v,w,x,y,z]$.
\end{remark}

\begin{proof}[Proof of Theorem~\ref{reg}]

Since $F_3G_3 \in (F_1,F_2)R$ and $F_1, F_2, F_3$ form a regular sequence, there
exist $\alpha, \beta \in R$ such that $G_3 = \alpha F_1 + \beta F_2$. In
\cite[\S 2, 3]{Si} we showed that in the polynomial ring $\mathbb{Z}[A,B,T]$
$$
\frac{(A+B)^q+(-A)^q+(-B)^q}{p} \Big[ (A+B)AB \Big]^k
$$
is an element of the ideal 
\begin{multline*}
(A+B)^{q+k}(T,A)^k (T,B)^k \ + \ A^{q+k}(T,B)^k (T+B, A+B)^k \\ 
+ \ B^{q+k}(T,A)^k (T-A, A+B)^k 
\end{multline*}
when $k=q-1$. We shall use this fact with $A=-F_2G_2$, \ $B=-F_3G_3$ \ and \ 
$T=\beta F_2F_3$. Note that this gives $A+B = F_1G_1, \ T+B = -\alpha F_1 F_3$, 
\ and \ $T-A = -F_1 G_1 - \alpha F_1 F_3$. Consequently
$$
\frac{(F_1G_1)^q+(F_2G_2)^q+(F_3G_3)^q}{p} \Big[ F_1G_1 F_2G_2 F_3G_3 \Big]^k 
$$
is an element of the ideal 
\begin{multline*}
(F_1G_1)^{q+k} (F_2F_3, \ F_2G_2)^k (F_2F_3, \ F_3G_3)^k \\
+ (F_2G_2)^{q+k} (F_2F_3, \ F_3G_3)^k (F_1 F_3, \ F_1G_1)^k \\ 
+ (F_3G_3)^{q+k} (F_2F_3, \ F_2G_2)^k (F_1 G_1 + \alpha F_1 F_3, \ F_1G_1)^k \\
\subseteq \ (F_1F_2F_3)^k(G_1^{q+k}, \ G_2^{q+k}, \ G_3^{q+k}).
\end{multline*}
The required result follows from this statement.
\end{proof}

\section{A Pl\"ucker relation}

One of the main goals of \cite{Si} was to settle Conjecture~\ref{conj} for the 
relation 
$$
u(vz-wy) + v(wx-uz) + w(uy-vx) = 0 \qquad \qquad (1)
$$
in the polynomial ring ${\mathbb Z}[u,v,w,x,y,z]$. This was accomplished by
establishing that for every prime integer $p$ and $q=p^e$ we have, for $k=q-1$, 
\begin{multline*}
\frac{1}{p} \left[ u^q(vz-wy)^q+v^q(wx-uz)^q+w^q(uy-vx)^q \right] \\
\times \left[ (vz-wy)(wx-uz)(uy-vx) \right]^k \\
 \in \ \Big( (vz-wy)^{q+k}, \ (wx-uz)^{q+k}, \ (uy-vx)^{q+k} \Big). 
\end{multline*}
The syzygy of the matrix 
$\displaystyle{\begin{pmatrix}vz-wy & wx-uz & uy-vx \end{pmatrix}}$ is
$$
\begin{pmatrix}
x & u\\
y & v\\ 
z & w\\ 
\end{pmatrix}
$$
but, by symmetry, Conjecture~\ref{conj} also holds for the relation 
$$
x(vz-wy) + y(wx-uz) + z(uy-vx) = 0. \qquad \qquad (2)
$$
In the light of Proposition~\ref{sum}, this gives us the following:

\begin{theorem}
Conjecture~\ref{conj} holds for 
$$
\Delta_1=vz-wy, \ \Delta_2=wx-uz, \ \Delta_3=uy-vx \ \in \
R={\mathbb Z}[u,v,w,x,y,z],
$$ 
i.e., if \ $F_1,F_2,F_3 \in R$ \ satisfy \ 
$F_1\Delta_1+F_2\Delta_2+F_3\Delta_3=0$ \
then for any prime power $q=p^e$ we have, for $k=q-1$,
$$
\frac{(F_1\Delta_1)^q + (F_2\Delta_2)^q + (F_3\Delta_3)^q}{p} 
(\Delta_1\Delta_2\Delta_3)^{k} 
\in (\Delta_1^{q+k}, \ \Delta_2^{q+k}, \ \Delta_3^{q+k})R.
$$
\end{theorem}

\bigskip

By taking a combination of $(1)$ and $(2)$ above, we get the relation
$$
(vz-wy)(ut-xs) + (wx-uz)(vt-ys) + (uy-vx)(wt-zs) = 0 
$$
in $R={\mathbb Z}[s,t,u,v,w,x,y,z]$. This is, of course, the Pl\"ucker relation 
$$
-\Delta_{34}\Delta_{12}+\Delta_{24}\Delta_{13}-\Delta_{23}\Delta_{14}
$$
where $\Delta_{ij}$ is the size two minor formed by picking rows $i$ and $j$ of
the matrix
$$
\begin{pmatrix} 
s & u & v & w \\ 
t & x & y & z \\ 
\end{pmatrix}.
$$
The syzygy of the matrix 
$\displaystyle{\begin{pmatrix}ut-xs & vt-ys & wt-zs \end{pmatrix}}$ is 
$$
\begin{pmatrix}
vz-wy & 0 & -wt+zs & vt-ys \\
wx-uz & -wt+zs & 0 & -ut+xs \\ 
uy-vx & vt-ys & ut-xs & 0
\end{pmatrix}.
$$
Since the Koszul relations are easily treated, the following theorem shows 
that Conjecture~\ref{conj} holds for $ut-xs, \ vt-ys, \ wt-zs$.

\begin{theorem}
\label{plucker}
Consider the Pl\"ucker relation
$$
(vz-wy)(ut-xs) + (wx-uz)(vt-ys) + (uy-vx)(wt-zs) = 0 
$$
in the polynomial ring $R = {\mathbb Z}[s,t,u,v,w,x,y,z]$. Then for $k=q-1$, 
we have
\begin{multline*}
\frac{1}{p} \left[ (vz-wy)^q(ut-xs)^q + (wx-uz)^q(vt-ys)^q + (uy-vx)^q(wt-zs)^q
\right] \\
\times [(ut-xs)(vt-ys)(wt-zs)]^k \\
\in \ \Big( (ut-xs)^{q+k}, \ (vt-ys)^{q+k}, \ (wt-zs)^{q+k} \Big)R.
\end{multline*}
\end{theorem}

Towards the proof of this theorem, we first record some identities with
binomial coefficients. These identities can be proved using Zeilberger's
algorithm (see \cite{PWS}) and the Maple package {\tt EKHAD}, but we include
proofs for the sake of completeness.

When the range of a summation is not specified, it is assumed to extend over
all integers. We set $\displaystyle{\binom{k}{i}=0}$ if $ i<0$ or if $k<i$.

\begin{lemma}
\label{binomiden}
\begin{alignat*}2
& \sum_{n} (-1)^n \binom{m+s-r}{m-n}\binom{k-n}{k-r} =
\begin{cases}
(-1)^{r-s} \binom{k-m}{s}, & \text{if } \ m \le k-s, \\
(-1)^r \binom{m+s-k-1}{s}, & \text{if } \ k+1 \le m, \\
0 & \text{else.}
\tag 1 
\end{cases} \\
& \sum_{r} (-1)^r \binom{2k-r}{m-1} \binom{m}{r-s} = 
0 \qquad \text{if} \qquad 1 \le m \le 2k-s. 
\tag 2 
\end{alignat*}
\end{lemma}

\begin{proof}
(1) Let 
$$
F(s,n) = (-1)^n \binom{m+s-r}{m-n}\binom{k-n}{k-r}, \qquad H(s)=\sum_n F(s,n),
$$
and $\displaystyle{G(s,n) = (-1)(k+1-n)F(s,n)}$. It is easily verified that
$$
\frac{G(s,n+1)}{G(s,n)} = \frac{(-1)(m-n)(r-n)}{(k+1-n)(s-r+n+1)}, 
\quad \frac{F(s+1,n)}{F(s,n)} = \frac{m+s+1-r}{n+s+1-r},
$$
and these can then be used to obtain the relation
$$
(s+1)F(s+1,n) - (m+s-k)F(s,n)= G(s,n+1) - G(s,n). 
$$
Summing the above equation with respect to $n$ gives
$$
(s+1)H(s+1) - (m+s-k)H(s) = 0,
$$
and using this recurrence relation for $H(s)$ we get
$$
H(s) = \frac{(m-k+s-1) \cdots (m-k)}{s \cdots 1} H(0).
$$
The result now follows since $H(0) = (-1)^r$.

\bigskip

(2) Let $\displaystyle{
F(r) = (-1)^r \binom{2k-r}{m-1} \binom{m}{r-s} 
}$
and $G(r) = (2k+1-r)(r-s) F(r)$. It is a routine verification that
$$
G(r) - G(r+1) = m(2k+1-m-s) F(r). 
$$
Consequently if $1 \le m \le 2k-s$ then
$$
\sum_{r} F(r) = \frac{1}{m(2k+1-m-s)} \sum_{r} \big(G(r) - G(r+1)\big) = 0.
$$
\end{proof}

We use the above results to establish the following equational identity:

\begin{lemma}
\label{asym}
\begin{multline*}
\sum_{r=0}^k \sum_{n=0}^r \binom{2k+1-r}{n}\binom{k-n}{r-n} 
\frac{T^{k-r}}{2k+1-r} (T-Y)^n (Y-Z)^{2k+1-r-n}Y^rZ^r \\
= \sum_{s=0}^k (-1)^s \binom{k}{s} \frac{T^{k-s}}{2k+1-s}
\big(Y^{2k+1}Z^s- Z^{2k+1}Y^s\big), 
\end{multline*}
\end{lemma}

\begin{proof}
Comparing the coefficients of $T^{k-s}$, we need to show that
\begin{multline*}
\sum_{r=s}^k \sum_{n=r-s}^r \frac{(-1)^{n-r+s}}{2k+1-r} 
\binom{2k+1-r}{n} \binom{k-n}{r-n} \binom{n}{r-s} \\
\times \ (Y-Z)^{2k+1-r-n}Y^{n+s} Z^r 
\ = \ \frac{(-1)^s}{2k+1-s} \binom{k}{s} \big(Y^{2k+1}Z^s - Z^{2k+1}Y^s\big), 
\end{multline*}
i.e., that
\begin{multline*}
\sum_{r=s}^k \sum_{n=r-s}^r \frac{(-1)^{n-r}}{2k+1-r} 
\binom{2k+1-r}{n} \binom{k-n}{r-n} \binom{n}{r-s} \\
\times \ (Y-Z)^{2k+1-r-n}Y^n Z^{r-s} \ 
= \ \frac{1}{2k+1-s} \binom{k}{s} \big( Y^{2k+1-s} - Z^{2k+1-s} \big).
\end{multline*}
The coefficient of $Y^mZ^{2k+1-s-m}$ in the expression on the left hand side is 
\begin{multline*}
\sum_{r,n} \frac{(-1)^{n+m+1}}{2k+1-r} \binom{2k+1-r-n}{m-n} \binom{2k+1-r}{n} 
\binom{k-n}{r-n} \binom{n}{r-s} \\
=\sum_r \frac{(-1)^{m+1}}{2k+1-r} \binom{m}{r-s} \binom{2k+1-r}{m} 
\sum_n (-1)^n \binom{m+s-r}{m-n} \binom{k-n}{r-n}. 
\end{multline*}
which we denote by $\gamma_{m,s}$. Note that $\gamma_{m,s}=0$ unless 
$m \le k-s$ or $m \ge k+1$ since 
$\displaystyle{
\sum_n (-1)^n \binom{m+s-r}{m-n} \binom{k-n}{r-n} = 0
}$ for such $m$ by Lemma~\ref{binomiden}~(1).

\bigskip

We next consider the case $m \le k-s$. It follows from Lemma~\ref{binomiden}~(1)
that 
$$
\gamma_{m,s} = \sum_r \frac{(-1)^{m+1+r-s}}{2k+1-r} \binom{m}{r-s} 
\binom{2k+1-r}{m} \binom{k-m}{s}.
$$
If $m=0$, then the only term that contributes to this sum is when $r=s$, and 
we get $\displaystyle{\gamma_{0,s}=\frac{-1}{2k+1-s}\binom{k}{s}}$. We may 
now assume $1 \le m \le k-s$ and then 
$$
\gamma_{m,s}= \frac{(-1)^{m+1-s}}{m} \binom{k-m}{s} 
\sum_r (-1)^r \binom{m}{r-s} \binom{2k-r}{m-1} =0 
$$
by Lemma~\ref{binomiden}~(2).

\bigskip

For $m \ge k+1$, Lemma~\ref{binomiden}~(1) gives 
\begin{align*}
\gamma_{m,s} & = \sum_r \frac{(-1)^{m+1+r}}{2k+1-r} \binom{m}{r-s} 
\binom{2k+1-r}{m} \binom{m+s-k-1}{s} \\
& = \frac{(-1)^{m+1}}{m} \binom{m+s-k-1}{s} 
\sum_r (-1)^r \binom{m}{r-s} \binom{2k-r}{m-1}.
\end{align*}
By Lemma~\ref{binomiden}~(2), this is zero if $m \le 2k-s$. If $m=2k+1-s$ the 
only term that contributes to this sum is when $r=s$, and
$\displaystyle{\gamma_{2k+1-s,s}=\frac{1}{2k+1-s}\binom{k}{s}}$. 
\end{proof}

This lemma enables us to establish the following crucial identity:

\begin{lemma}
\label{iden}
\begin{align*}
\sum_{r=0}^k \sum_{n=0}^r \frac{k+1}{2k+1-r} & \binom{2k+1-r}{n} \binom{k-n}{r-n} 
\ T^{k-r} \\ 
\times \bigg[ & X^{2k+1}(T-Y)^n (Y-Z)^{2k+1-r-n}Y^rZ^r \\
+ & Y^{2k+1}(T-Z)^n (Z-X)^{2k+1-r-n}Z^rX^r \\
+ & Z^{2k+1}(T-X)^n (X-Y)^{2k+1-r-n}X^rY^r \bigg] = 0
\end{align*}
\end{lemma}

\begin{proof}
Using Lemma~\ref{asym}, the left hand side in the equation equals
\begin{multline*}
\sum_{s=0}^k (-1)^s \binom{k}{s} \frac{k+1}{2k+1-s}T^{k-s} 
\times \bigg[
X^{2k+1}(Y^{2k+1}Z^s- Z^{2k+1}Y^s) \\ \ + \ 
Y^{2k+1}(Z^{2k+1}X^s- X^{2k+1}Z^s) \ + \
Z^{2k+1}(X^{2k+1}Y^s- Y^{2k+1}X^s)
\bigg] \ = \ 0.
\end{multline*}
\end{proof}

\begin{proof}[Proof of Theorem~\ref{plucker}]
Since $k=q-1$, we have $(k+1)/(2k+1-r)=q/(2q-1-r)$. For $0 \le r \le q-2$, if 
$q/(2q-1-r) = a/b$ for relatively prime integers $a$ and $b$, then $p$ divides 
$a$ since $q+1 \le 2q-1-r \le 2q-1$. Consequently there exists an integer
$d \in {\mathbb Z}$ such that $d$ is relatively prime to $p$ and
$$
\frac{dq}{2q-1-r} \in {\mathbb Z} \quad \text{ for all } \quad 0 \le r \le q-1.
$$
After replacing $d$ by a suitable multiple, if necessary, we may 
assume that $d \equiv 1 \mod p$. Using Lemma~\ref{iden} with 
$$
T=\frac{t}{s}, \qquad X=\frac{t}{s}-\frac{x}{u}, 
\qquad Y=\frac{t}{s}-\frac{y}{v}, \qquad Z=\frac{t}{s}-\frac{z}{w},
$$
we get
\begin{align*}
\sum_{r=0}^k & \sum_{n=0}^r \frac{k+1}{2k+1-r} \binom{2k+1-r}{n} 
\binom{k-n}{r-n} \ \left( \frac{t}{s} \right)^{k-r} \\ 
\times \bigg[ 
& \left( \frac{ut-xs}{su} \right)^{2k+1} \left(\frac{y}{v} \right)^n 
 \left( \frac{vz-wy}{vw} \right)^{2k+1-r-n} \left(\frac{vt-ys}{sv} \right)^r 
 \left( \frac{wt-zs}{sw} \right)^r \\
& \left( \frac{vt-ys}{sv} \right)^{2k+1} \left(\frac{z}{w} \right)^n 
 \left( \frac{wx-uz}{wu} \right)^{2k+1-r-n} \left(\frac{wt-zs}{sw} \right)^r 
 \left( \frac{ut-xs}{su} \right)^r \\
& \left( \frac{wt-zs}{sw} \right)^{2k+1} \left(\frac{x}{u} \right)^n 
 \left( \frac{uy-vx}{uv} \right)^{2k+1-r-n} \left(\frac{ut-xs}{su} \right)^r 
 \left( \frac{vt-ys}{sv} \right)^r \bigg] \\
& = 0.
\end{align*}
Multiplying by $ds^{4k+1}(uvw)^{2k+1}$ clears denominators, and gives us the
equational identity (with integer coefficients):
\begin{multline*}
\sum_{r=0}^k \sum_{n=0}^r \frac{d(k+1)}{2k+1-r} \binom{2k+1-r}{n} 
\binom{k-n}{r-n} \ (st)^{k-r} \\ \times \bigg[ 
(ut-xs)^{2k+1} (wy)^n (vz-wy)^{2k+1-r-n} (vt-ys)^r (wt-zs)^r \\ + \
(vt-ys)^{2k+1} (uz)^n (wx-uz)^{2k+1-r-n} (wt-zs)^r (ut-xs)^r \\ + \
(wt-zs)^{2k+1} (vx)^n (uy-vx)^{2k+1-r-n} (ut-xs)^r (vt-ys)^r 
\bigg] = 0.
\end{multline*}
Note that for $0 \le r \le k=q-1$ and $0 \le n \le r$,
$$
\frac{d(k+1)}{2k+1-r} \binom{2k+1-r}{n} \binom{k-n}{r-n} \equiv
\begin{cases}
1 \mod p, & \text{if } \ r=k \text{ and } n=0, \\
0 \mod p, & \text{else}.
\end{cases}
$$
Consequently we get
\begin{align*}
 & (ut-xs)^{2k+1} (vz-wy)^{k+1} (vt-ys)^k (wt-zs)^k \\ 
+ \ & (vt-ys)^{2k+1} (wx-uz)^{k+1} (wt-zs)^k (ut-xs)^k \\ 
+ \ & (wt-zs)^{2k+1} (uy-vx)^{k+1} (ut-xs)^k (vt-ys)^k \\
& \qquad \qquad \in \ \Big( p(ut-xs)^{2k+1}, \ p(vt-ys)^{2k+1}, \ 
p(wt-zs)^{2k+1} \Big) R.
\end{align*}
Finally, this shows that
\begin{multline*}
\frac{1}{p}[(vz-wy)^q(ut-xs)^q + (wx-uz)^q(vt-ys)^q + (uy-vx)^q(wt-zs)^q] \\
\times [(ut-xs)(vt-ys)(wt-zs)]^k \\
\in \ \Big( (ut-xs)^{q+k}, \ (vt-ys)^{q+k}, \ (wt-zs)^{q+k} \Big)R.
\end{multline*}
\end{proof}

\section*{Acknowledgments} 

I am indebted to Mel Hochster for valuable discussions regarding this material.
The use of the Maple package {\tt EKHAD} is gratefully acknowledged. I also
want to express my appreciation to Gennady Lyubeznik and Xavier G\'omez-Mont
Avalos for organizing the {\it Conference on Local Cohomology}.

\end{document}